\documentclass{amsart}

\newtheorem{teo}{Theorem}[section]

\newtheorem{corol}[teo]{Corollary}

\theoremstyle{remark}
\newtheorem{rem}{Remark}[section]

\newcommand{\nc}{\newcommand}
\nc{\sig}{\sigma}
\nc{\om}{\omega}
\nc{\ti}{\theta}

\newcommand{\al}{\alpha}

\newcommand{\hc}{$\{J_{\alpha}\}_{\alpha =1,2,3}$}
\nc{\hk}{\tilde{g}}
\nc{\ex}{\wedge}

\nc{\cc}{\Bbb C}
\nc{\hh}{\Bbb H}
\nc{\pp}{\Bbb P}
\nc{\rr}{\Bbb R}

\newcommand{\cg}{{\frak g}}
\nc{\ch}{{\mathcal H}}

\nc{\nag}{\nabla\sp g}
\nc{\nah}{\nabla\sp{\ch}} 
\numberwithin{equation}{section}

\title[Hyper-Hermitian metrics 
which are conformally  hyper-K\"ahler]{Homogeneous hyper-Hermitian metrics 
which are conformally  hyper-K\"ahler} 

\author[MARIA LAURA BARBERIS]{Mar\'{\i}a Laura Barberis}

\address{FaMAF, Universidad Nacional de C\'ordoba, 
Ciudad Universitaria, 5000 - C\'ordoba, Argentina}

\email{barberis@mate.uncor.edu}

\thanks{The author was partially supported by CONICET, ESI (Viena) and FOMEC (Argentina)}



\begin{document}

\maketitle

\begin{abstract}
 Let $g$ be a  hyper-Hermitian metric on a simply connected 
hypercomplex four-manifold $(M, \ch)$. 
We show that when the isometry group $I(M,g)$ contains 
a subgroup acting simply transitively on $M$ by hypercomplex isometries   
 then the metric $g$ is conformal to a hyper-K\"ahler metric. 
We describe explicitely the corresponding hyper-K\"ahler metrics and 
it follows that, in four dimensions,  these are the only hyper-K\"ahler metrics containing a homogeneous metric in its conformal class.    
\end{abstract}

\section{Preliminaries}
A  hypercomplex structure on a $4n$-dimensional 
manifold $M$ is a family $\ch=$ \hc\  
of fibrewise endomorphisms of the tangent bundle $TM$ of $M$ 
 satisfying:
\begin{gather}
\label{Quat}
J_{\al}^2=-I, \hspace{.5cm}  \al=1,2,3 \hspace{1cm}
J_{1}J_{2}=-J_{2}J_{1}=J_3, \hspace{.5cm}     \\
   \label{intr} N_{\al}\equiv 0, \hspace{.8cm}\al=1,2,3 
\end{gather}
 where
$I$\ is the identity on the tangent space
$T_pM$\  of $M$ at $p$ for all $p$ in $M$ and $N_{\al}$\ is the 
Nijenhuis tensor
corresponding to $J_{\al}$: \[
N_{\al}(X,Y)=[J_{\al}X,J_{\al}Y]-[X,Y] 
 - J_{\al}([X,J_{\al}Y]+[J_{\al}X,Y])  \]
for all $\,X,Y\,$ vector fields on $M$. A differentiable map 
$f:M \rightarrow M$ is said to be hypercomplex if it is holomorphic 
with respect to $J_{\al}, \, \al=1,2,3$.  The group of hypercomplex diffeomorphisms on ($M, \ch)$ will be denoted by  Aut$(\ch)$. 

A riemannian metric $g$ on a 
hypercomplex manifold $(M, \ch)$ is called hyper-Hermitian when 
$g(J_{\al}X,J_{\al}Y)=g(X,Y)$ for all vectors fields 
$X,Y$ on $M$, $ \al=1,2,3$. 

Given a  manifold $M$ with a hypercomplex structure $\ch=$ \hc\ and a 
hyper-Hermitian metric $g$ consider  the $2$-forms $\om_{\al}, 
\al=1,2,3$, 
defined by 
\begin{equation} \label{2forms}
\om_{\al}(X,Y)=g(X,J_{\al}Y). \end{equation}
The metric $g$ is said to be hyper-K\"ahler when $d\om_{\al}=0$ 
for $\al=1,2,3$.  

It is well known 
that a hyper-Hermitian metric $g$ 
is conformal to a hyper-K\"ahler metric $\hk$ if and only if 
there exists an exact $1$-form $\ti\in \Lambda ^1M$ such that 
\begin{equation}  \label{lck}
d\om_{\al}=\ti \ex \om_{\al},  \hspace{,5cm} \al=1,2,3
\end{equation} 
where, if $g= e^f \, \hk$ for some $f \in C^{\infty}(M)$, then $\ti=df$. 

We prove the following result:
\begin{teo} \label{Thm}
Let $(M, \ch, g)$ be a simply connected   hyper-Hermitian 
$4$-manifold. Assume that there exists a Lie group $G \subset I(M,g) \, 
\cap \, \text{Aut} \, ( \ch)$ 
acting simply transitively on $M$.  
Then $g$ is conformally hyper-K\"ahler.  
\end{teo}
We conclude that one of the hyper-K\"ahler metrics constructed by the Gibbons-Hawking ansatz \cite{gib} contains a homogeneous hyper-Hermitian metric 
in its conformal class. This hyper-Hermitian metric is not symmetric and has negative sectional curvature \cite{bar}. 

As a  consequence of  Theorem \ref{Thm} and the results in \cite{bar} we obtain that the following 
symmetric riemannian metrics are conformally hyper-K\"ahler:
\begin{itemize}
\item the riemannian product of the canonical metrics on ${\Bbb R} \times S^3$;  
\item the riemannian product of the canonical metrics on 
  ${\Bbb  R} \times {\Bbb R} H^3 $, where  ${\Bbb R} H^3$ denotes the 
real hyperbolic space;
\item the canonical metric on the real hyperbolic space 
 ${\Bbb R} H^4$.
\end{itemize}

\begin{small} {\bf Acknowledgements.} I would like to thank the organizers of the program 
{\it Holonomy Groups in Differential Geometry} for giving me the opportuniy 
to visit the {\it Erwin Schr\"odinger Institute}, Vienna.   
 I am also grateful to D. 
Alekseevsky, I. Dotti, L. Ornea and S. Salamon for useful conversations. 
\end{small}

\section{Proof of the main theorem} \label{sec2}

{\it Proof of Theorem \ref{Thm}.}  
Since $G$ acts simply transitively on $M$ then $M$ is diffeomorphic to 
$G$ and therefore the hypercomplex structure  and hyper-Hermitian metric 
can be transferred to $G$ and  will also be denoted by 
\hc\ and $g$, respectively. Since $G$ acts 
by hypercomplex isometries it follows that both \hc\ and $g$ are left 
invariant on $G$. All such simply connected Lie groups were classified in 
\cite{bar}, where it is shown that the Lie algebra $\cg$ of $G$ is either 
abelian or isomorphic 
to one of the following Lie algebras (we fix  an 
orthonormal basis $\{ e_j\}_{j=1,\dots ,4}$
of $\cg$):
\begin{enumerate}
\item  $[e_2,e_3]=e_4,\;\;\;  [e_3,e_4]=e_2,\;\;\; [e_4,e_2]=e_3,\;\;\;e_1$ central;
\item $[e_1,e_3]=e_1, \;\; [e_2,e_3]=e_2, \;\; [e_1,e_4]=e_2, 
\;\; [e_2,e_4]=-e_1$;
\item $[e_1, e_j]=e_j, \; j=2,3,4$;
\item $[e_3,e_4]= \frac 12 \, e_2, \;\; [e_1,e_2]=e_2, \;\; 
[e_1,e_j]= \frac 12 \, e_j, \; j=3,4$. 
\end{enumerate}
Observe that in case 1 above $M$ is diffeomorphic to ${\Bbb R} \times S^3$ while 
in the remaining cases it is diffeomorphic to  ${\Bbb R}^4$, therefore in all cases any closed form on $M$ is exact.  
We now proceed by finding in each case a closed form $\ti \in \Lambda^1\cg ^*$ 
satisfying  \eqref{lck}. Note that we work on the Lie algebra level since  
$g$ and $\om_{\al}$ are all left invariant on $G$.  
Let 
 $\{ e^j\}_{j=1,\dots ,4} \subset \Lambda^1\cg^*$ 
be the dual basis of  $\{ e_j\}_{j=1,\dots ,4}$. From now on 
we will write 
$e^{ij \cdots}$  to denote $ e^i \ex e^j \ex \cdots$.   In all the cases below the $2$-forms $\om_{\al}$ are determined from \eqref{2forms} in terms of the hypercomplex structures constructed in \cite{bar}. 

\

Case 1. The $2$-forms $\om_{\al}$ are given as follows:
\[ \om _1= -e^{12}-e^{34}, \hspace{.8cm}  
\om _2= -e^{13}+e^{24}, \hspace{.8cm} 
\om _3= -e^{14}-e^{23}.  
\]
To calculate $d\om_{\al}$ we obtain first $de^j$ (recall that 
$d\sigma (x,y)=-\sigma[x,y]$ for $\sigma \in \Lambda^1\cg ^*$):
\begin{equation} \label{eq1}
de^1=0, \hspace{.6cm} de^2= -e^{34}, \hspace{.6cm} 
 de^3= e^{24}, \hspace{.6cm}
 de^4= -e^{23}. \end{equation}
These equations and  the fact that $d(\sigma \ex \tau)=
d\sigma \ex \tau + (-1)^r \sigma \ex d\tau$ for all $\sigma \in 
\Lambda ^r \cg^*$ give the following formulas:
\[ 
 d\om_1= -e^{134}, \hspace{.8cm} 
 d\om_2= e^{124}, \hspace{.8cm} 
 d\om_3= -e^{123}  
\]
from which we conlude that \eqref{lck} holds  for 
$\ti= e^1$, which is closed and therefore exact since 
$G$ is diffeomorphic to ${\Bbb R} \times S^3$. We conclude that this hyper-Hermitian metric, 
which, as shown in \cite{bar}, is homothetic to the riemannian 
product of the canonical metrics  
on ${\Bbb R} \times S^3$,  
 is conformal to a hyper-K\"ahler metric. 

\

Case 2. In this case we have the following equations for  
$\om_{\al}$:
\[ \om _1= e^{14}-e^{23}, \hspace{.8cm}  
\om _2= -e^{12}+e^{34}, \hspace{.8cm} 
\om _3= -e^{13}-e^{24}.  
\]
and we calculate
\begin{gather} \label{eqaff}
de^1=-e^{13}+e^{24}, \hspace{.5cm} de^2= -e^{23}-e^{14}, \hspace{.5cm} 
 de^3= 0, \hspace{.6cm}
 de^4= 0, \\ 
 d\om_1= -2e^{134}, \hspace{.8cm} 
 d\om_2= -2e^{123}, \hspace{.8cm} 
 d\om_3= 2e^{234}  
\end{gather} 
so that \eqref{lck} is satisfied for $\ti=2e^3$, which again is closed, so this 
hyper-Hermitian metric is also conformal to a hyper-K\"ahler metric. 
In this case the hyper-Hermitian metric is homothetic to the 
riemannian product of the canonical metrics on 
  ${\Bbb  R} \times {\Bbb R} H^3 $, where  ${\Bbb R} H^3$ denotes the 
real hyperbolic space. 

\

Case 3. 
In this case the $2$-forms  
$\om_{\al}$ are given as follows:
\[ \om _1= -e^{12}-e^{34}, \hspace{.8cm}  
\om _2= -e^{13}+e^{24}, \hspace{.8cm} 
\om _3= -e^{14}-e^{23}  
\]
and a calculation of exterior derivatives gives:
 \begin{gather}
de^1=0, \hspace{1cm} de^j= -e^{1j}, \hspace{.3cm} j=2,3,4 \\
 d\om_1= 2e^{134}, \hspace{.8cm} 
 d\om_2= -2e^{124}, \hspace{.8cm} 
 d\om_3= -2e^{123}  
\end{gather} 
so that \eqref{lck} is satisfied for $\ti=-2e^1$. 
This hyper-Hermitian metric is homothetic to the canonical 
metric on the real hyperbolic space ${\Bbb R}H^4$. 

\

Case 4.  In this case we have the following equations for  
$\om_{\al}$:
\[ \om _1= -e^{12}+e^{34}, \hspace{.8cm}  
\om _2= -e^{13}-e^{24}, \hspace{.8cm} 
\om _3= e^{14}-e^{23}  
\]
and we calculate
 \begin{gather}\label{eqch2}
de^1=0, \hspace{.5cm} 
de^2=-e^{12} -\frac 12 \, e^{34},   \hspace{.5cm}
de^j= - \frac 12 \, e^{1j}, \hspace{.3cm} j=3,4 \\
 d\om_1= -\frac 32 \, e^{134}, \hspace{.8cm} 
 d\om_2= \frac 32 \, e^{124}, \hspace{.8cm} 
 d\om_3= \frac 32 \, e^{123}  
\end{gather} 
so that \eqref{lck} is satisfied for $\ti=-\frac 32\, e^1$. 
This hyper-Hermitian metric is not symmetric and has negative sectional curvature (cf. \cite{bar}).

\begin{rem}
All the hyper-Hermitian manifolds $(M,\ch , g)$ considered above 
admit a connection $\nabla$ such that: 
\[
\nabla g=0, \hspace{1cm} \nabla J_{\al}=0, \hspace{.3cm} \al=1,2,3 
\]
and the $(3,0)$ tensor $c(X,Y,Z)=g(X, T(Y,Z))$ is totally skew-symmetric,
 where $T$ is the torsion of $\nabla$. Such a connection is called an 
HKT connection (cf. \cite{poon}). In case $M$ is diffeomorphic to 
 ${\Bbb R} \times S^3$ it can be shown that, moreover, the 
corresponding $3$-form $c$ is closed.  
\end{rem} 

\section{Coordinate description of the Hyper-K\"ahler metrics}
In this section we will use global coordinates on each of the Lie groups considered in the previous section to describe the corresponding hyper-K\"ahler metrics. This will allow us to 
identify the hyper-K\"ahler metric in \S \ref{sec2}, Case 4, with one constructed by the Gibbons-Hawking ansatz \cite{gib}. 

\

Case 1.   $G= \hh ^*=GL(1,\hh)=\left\{ \;
\left( \begin{matrix}
x & -y& -z& -t \\
y & x& -t& z\\
z& t& x& -y\\
t& -z& y& x
\end{matrix} \right) \; :\; (x,y,z,t) \in \rr^4 - \{0\} \; 
\right\} $.  We obtain a basis of left invariant $1$-forms on $G$ as follows. 
Set $r^2=x^2+y^2+z^2+t^2$, $r>0$, and  $\Omega= g^{-1}dg$ for $g\in G$, that is, \[
\text{if }   \;\; g= \left( \begin{matrix}
x & -y& -z& -t \\
y & x& -t& z\\
z& t& x& -y\\
t& -z& y& x
\end{matrix} \right) \hspace{.3cm}  \text{then }\hspace{.3cm}
\Omega=  \left( \begin{matrix}
\sig_1 & -\sig_2& -\sig_3& -\sig_4 \\
 \sig_2& \sig_1& -\sig_4& \sig_3\\
\sig_3& \sig_4& \sig_1& -\sig_2\\
\sig_4& -\sig_3& \sig_2& \sig_1
\end{matrix} \right)
\]
where 
\[ \left( \begin{matrix}
\sig_1 \\ \sig_2\\ \sig_3\\ \sig_4\end{matrix} \right)= 
\frac{1}{r^2} \left( \begin{matrix}
x & y& z& t \\
-y & x& t& -z\\
-z& -t& x& y\\
-t& z& -y& x
\end{matrix} \right)\left( \begin{matrix}
dx \\ dy\\ dz\\ dt \end{matrix} \right). 
\]
Then $\sig_j, \; 1\leq j\leq 4$, is a basis of left invariant $1$-forms on $G$ and it follows from $d\Omega + \Omega \ex \Omega=0$ that 
\[
d\sig _1=0, \hspace{.6cm} d\sig_2= -2\sig_3\ex\sig_4, \hspace{.6cm} 
 d\sig_3= 2\sig_2\ex \sig_4, \hspace{.6cm}
 d\sig_4= -2\sig_2 \ex\sig_3. \]
Setting \[ e^1=2\sig_1,\hspace{.5cm} e^2=2\sig_2,\hspace{.5cm}  e^3=2\sig_3,
\hspace{.5cm} e^4=2\sig_4, \]
 so that $\{e^j\}_{1\leq j\leq 4}$ satisfy 
\eqref{eq1},   
 the left-invariant hyper-Hermitian metric is 
\begin{equation} 
g=(e^1)^2+(e^2)^2+(e^3)^2+(e^4)^2= \frac 4{r^2}(dx^2+dy^2+dz^2+dt^2)
\end{equation} 
that is, $g$ is the standard conformally flat 
metric on $\rr^4 - \{0\}$,  
and since the Lee form is $\ti=e^1= d(2\log r)$ the corresponding hyper-K\"ahler metric is $\hk=e^{-2\log r}g$, that is, 
\begin{equation} 
\hk= \frac 4{r^2}\left( \frac{(dr)^2}{r^2}+ (\sig_2)^2
+(\sig_3)^2+(\sig_4)^2\right)=\frac 4{r^4}(dx^2+dy^2+dz^2+dt^2). 
\end{equation}  
Observe that the standard metric on any coordinate quaternionic Hopf  surface is locally conformally equivalent to $\hk$ (cf. \cite{boy}).   

\

Case 2. Define a product on $\rr ^4$ as follows:
\begin{equation*} \begin{split}
(x,y,z,t)&(x',y',z',t')= \\&(x+e^z(x'\cos t-y'\sin  t), y+e^z(x'\sin t+y'\cos t), z+z',t+t').\end{split}
\end{equation*}
This defines a Lie group structure on $\rr ^4$  that makes it isomorphic to the Lie group considered in \S \ref{sec2}, Case 2.  The following $1$-forms are left-invariant with respect to the above    product:
\begin{equation*} \begin{align}   e^1&=e^{-z}\cos t \, dx +e^{-z} \sin t \, dy, 
&e^3&=-dz,\\ 
   e^2&=-e^{-z} \sin t \, dx+e^{-z}\cos t \, dy, & 
  e^4&=-dt
\end{align}\end{equation*}
These forms satisfy  relations \eqref{eqaff}. The hyper-Hermitian metric is therefore given as follows:
\[   g= (e^1)^2+(e^2)^2+(e^3)^2+(e^4)^2= e^{-2z} (dx^2+dy^2) + dz^2+dt^2
\]
and the Lee form is $\ti= 2e^3=-2dz$, so that the hyper-K\"ahler metric becomes \[ 
\hk=e^{2z}g=(dx^2+dy^2)+ e^{2z}(dz^2+dt^2). \]
Observe that  the change of coordinates $s= e^z$ gives the following simple form for $\hk$ on $\rr^+ \times \rr^3$:
\[  \hk=  dx^2+dy^2+(ds^2+s^2dt^2).  \]
This allows us to identitify $\hk$ with the 
riemannian product of two K\"ahler metrics: the euclidean metric on $\rr^2$ with the 
warped product cone metric on $\rr^+\times \rr$ (cf. \cite{BG}). 

\

Case 3. We endow $\rr ^4$ with the following product:
\[ (x,y,z,t)(x',y',z',t')=(x+ e^t x', y+e^t y', z+ e^t z', t+t') \]
thereby obtaining the Lie group structure considered in \S \ref{sec2}, Case 3, with corresponding left-invariant $1$-forms:
\[ e^1= dt, \hspace{.5cm}
e^2=e^{-t}dx, \hspace{.5cm}
e^3=e^{-t}dy,\hspace{.5cm}
 e^4=e^{-t}dz.
\]
The hyper-Hermitian metric is therefore
\[   g= (e^1)^2+(e^2)^2+(e^3)^2+(e^4)^2= 
            e^{-2t}(dx^2+dy^2+dz^2) + dt^2
\]
with corresponding Lee form $\ti=-2e^1=-2dt$, yielding the following 
hyper-K\"ahler metric: 
\[     \hk= e^{2t}g= dx^2+dy^2+dz^2 + e^{2t}dt^2.      
\]
Setting $s=e^t$, $\hk$ is the euclidean metric  $ds^2+dx^2+dy^2+dz^2$  on $\rr^+ \times \rr^3$. 

\

Case 4. Consider the following product on  $\rr ^4$:
\[     (x,y,z,t)(x',y',z',t')= (x+e^{\frac t2} x', y+e^{\frac t2}y',
z+e^tz'+\frac{e^{\frac t2}}{4}(xy'-yx'), t+t')
\]
which yields the Lie group structure considered in \S \ref{sec2}, Case 4. It is easily checked that the following left-invariant $1$-forms satisfy \eqref{eqch2}:
\[  e^1=dt, \hspace{.4cm}e^2= e^{-t}(dz-\frac 14xdy+\frac 14ydx), \hspace{.4cm} e^3=e^{-\frac t2}dx, \hspace{.4cm} e^4=e^{-\frac t2}dy.  
\]
The hyper-Hermitian metric is now obtained as in the above cases:
 \begin{equation*}   \begin{split} g=(e^1)^2+(e^2)^2+&(e^3)^2+(e^4)^2\\
          =&dt^2+e^{-t}(dx^2 +dy^2) +e^{-2t} (dz-\frac 14(xdy-ydx))^2        
\end{split}\end{equation*}
and the Lee form is $\ti= -\frac 32 dt$, from which we obtain the hyper-K\"ahler metric as usual:
\[ \hk= e^{-\frac 32 t}dt^2 + e^{-\frac t2 } (dx^2 +dy^2) + 
e^{-\frac t2 }(dz-\frac 14(xdy-ydx))^2.   \]                             Setting $s=e^{\frac t2}$, $\hk$ becomes 
\[       \hk= s(ds^2+dx^2 +dy^2) + \frac 1s (dz-\frac 14(xdy-ydx))^2                             \]
on $\rr^+\times \rr^3$, which allows us to identify $\hk$ with 
one of the hyper-K\"ahler metrics constructed by the 
Gibbons-Hawking ansatz \cite{gib}. The identification 
is easily  obtained from  \cite{leb}, Proposition 1.

We can now rephrase Theorem \ref{Thm} as follows, where $[h]$ denotes the 
conformal class of $h$:

\begin{corol} Let $h$ be a hyper-K\"ahler metric on a simply connected 
hypercomplex $4$-manifold $(M, \ch)$ such that there exist $g\in [h]$ and a Lie group 
$G \subset I(M,g) \cap \text{Aut}\,( \ch)$ acting simply transitively on $M$. 
 Then ($M, h$) is homothetic to either $\rr^4$ with the euclidean metric or one of the following riemannian manifolds: 
\begin{enumerate}
\item $M= \rr^4 - \{ 0\}, \;\; h= r^{-4}(dx^2+dy^2+dz^2+dt^2),
$
\item $M=\rr^ 2 \times (\rr^+ \times \rr),\;\; h= (dx^2+dy^2) + (ds^2+ s^2 dt^2), $
\item $M=\rr^+ \times \rr^3
, \;\;h=ds^2+dx^2+dy^2+dz^2$, 
\item $M=\rr^+\times \rr^3, \;\;h= s(ds^2+dx^2 +dy^2) + s^{-1} (dz-\frac 14(xdy-ydx))^2$.
\end{enumerate}
\end{corol}

\bibliographystyle{amsplain}

\end{document}